\documentclass[conference]{IEEEtran}
\IEEEoverridecommandlockouts
\usepackage{cite}
\usepackage{amsmath,amssymb,amsfonts}
\usepackage{graphicx}
\usepackage{textcomp}
\usepackage{xcolor}
\usepackage{authblk}
\usepackage[linesnumbered,ruled,vlined]{algorithm2e}

\def\BibTeX{{\rm B\kern-.05em{\sc i\kern-.025em b}\kern-.08em
    T\kern-.1667em\lower.7ex\hbox{E}\kern-.125emX}}
\begin{document}

\title{Energy Efficient Altitude Optimization\\
 of an Aerial Access Point  \\
\thanks{©2020 IEEE. Personal use of this material is permitted. Permission from IEEE must be obtained for all other uses, in any
current or future media, including reprinting/republishing this material for advertising or promotional purposes, creating new
collective works, for resale or redistribution to servers or lists, or reuse of any copyrighted component of this work in other
works.
This version of the work has been accepted for publication in the IEEE 31st PIMRC’20 - Workshop on UAV Communications for 5G and Beyond. This work is supported by the project PAINLESS which has received funding from the European Union’s Horizon 2020 research and innovation programme under grant agreement No 812991.}
}

\author[1,2]{Nithin Babu}
\author[3]{Konstantinos  Ntougias}
\author[1,2]{Constantinos B Papadias}
\author[2]{Petar Popovski}
%\author[1]{Stavros Kotsopoulos}

\affil[1]{Research, Technology and Innovation Network (RTIN), 
The American College of Greece, Greece }
\affil[ ]{\textit {\{nbabu,cpapadias\}@acg.edu}}
\affil[2]{Department of Electronic Systems, Aalborg University, Denmark}
\affil[ ]{\textit {\{niba,cop,petarp\}@es.aau.dk}}
\affil[3]{University of Cyprus}
\affil[ ]{\textit {\{ntougias.konstantinos\}@ucy.ac.cy}}
\renewcommand\Authands{ and }

%\author{
%%\IEEEauthorblockN{Nithin Babu}
%%\IEEEauthorblockA{\textit{B-WiSe Lab, Athens Information Technology(AIT), Greece} \\
%%\textit{Dept. of Electronic Systems, Aalborg University, Denmark}\\
%%%Athens, Greece \\
%%niba@es.aau.dk, nithinbabu29@gmail.com}
%\IEEEauthorblockN{Constantinos B Papadias}
%\IEEEauthorblockA{\textit{B-WiSe Lab, Athens Information Technology(AIT), Greece } \\
%\textit{Adjunct Professor, Aalborg University, Denmark}\\
%  cop@es.aau.dk}
%\and
%\IEEEauthorblockN{Konstantinos  Ntougias}
%\IEEEauthorblockA{\textit{B-WiSe Lab, Athens Information Technology(AIT), Greece} \\
%kontou@ait.gr}
%\and
%\IEEEauthorblockN{Constantinos B Papadias}
%\IEEEauthorblockA{\textit{B-WiSe Lab, Athens Information Technology(AIT), Greece } \\
%\textit{Adjunct Professor, Aalborg University, Denmark}\\
%  cop@es.aau.dk}
%\and
%\IEEEauthorblockN{Petar Popovski}
%\IEEEauthorblockA{\textit{Professor, Aalborg University, Denmark}
%\\
%%City, Country \\
%petarp@es.aau.dk}
%%\and
%%\IEEEauthorblockN{5\textsuperscript{th} Given Name Surname}
%%\IEEEauthorblockA{\textit{dept. name of organization (of Aff.)} \\
%%\textit{name of organization (of Aff.)}\\
%%City, Country \\
%%email address}
%%\and
%%\IEEEauthorblockN{6\textsuperscript{th} Given Name Surname}
%%\IEEEauthorblockA{\textit{dept. name of organization (of Aff.)} \\
%%\textit{name of organization (of Aff.)}\\
%%City, Country \\
%%email address}
%}

\maketitle

\begin{abstract}
In this paper, we propose an energy-efficient optimal altitude for an aerial access point (AAP), which acts as a flying base station to serve a set of ground user equipment (UE). Since the ratio of total energy consumed by the aerial vehicle to the communication energy is very large, we include the aerial vehicle's energy consumption in the problem formulation. After considering the energy consumption model of the aerial vehicle, our objective is translated into a non-convex optimization problem of maximizing the global energy efficiency (GEE) of the aerial communication system, subject to altitude and minimum individual data rate constraints. At first, the non-convex fractional objective function is solved by using sequential convex programming (SCP) optimization technique. To compare the result of SCP with the global optimum of the problem, we reformulate the initial problem as a monotonic fractional optimization problem (MFP) and solve it using the polyblock outer approximation (PA) algorithm. Numerical results show that the candidate solution obtained from SCP is the same as the global optimum found using the monotonic fractional programming technique. Furthermore, the impact of the aerial vehicle's energy consumption on the optimal altitude determination is also studied.      
\end{abstract}
\section{Introduction}
The role of uninhabited AAP in the deployment of emergency networks such as deploying aerial base stations to provide reliable connectivity in disaster areas \cite{DISASTER} or in social events such as concerts is vital. In Japan, earthquake affected areas were provided with internet access with the help of unmanned aerial vehicles (UAV) \cite{japan}. Cellular coverage extension using drone deployed base stations by Nokia's F-cell technology is another proven application of portable access points \cite{nokia}. The mobility and ability of aerial vehicles to adjust their altitude to improve the probability of line-of-sight (LoS) communication channel to the ground UEs makes them suitable for acting as relays in the internet of things (IoT) applications \cite{IOT}. Despite all these applications, the efficiency of an aerial communication system (ACS) is highly dependent on the limited energy available at the aerial vehicle \cite{rui}. Any improvement in the energy efficiency of ACS  implies longer aerial vehicle hovering, hence more information bits transmitted to UEs.

Compared to the conventional cellular communication systems, the total energy required by ACS is very high. This is because, in ACS, in addition to the communication-related energy, the aerial vehicle consumes energy during vertical climb and hovering. Most of the works in the literature only consider communication-related energy, which is suboptimal in the case of an ACS. In \cite{angleforlos}, the authors present an analytical approach to optimize the altitude of low altitude aerial platforms to maximize the radio coverage area. The authors of \cite{Joint} jointly optimize the flying altitude and the antenna beamwidth for throughput maximization. A new 3-dimensional deployment plan for the drone-base station to serve the users based on their service requirements, while minimizing the number of drones, is presented in \cite{kalantari2016number}. The work in  \cite{lyu2016placement} proposes a new polynomial-time complex spiral mobile base station placement algorithm in UAV-UE communications. The works in \cite{alzenad20173},\cite{mozaffari2015drone} find the optimal altitude for UAV-base stations that maximizes the number of covered users using the minimum transmit power.

None of the above works consider the energy consumption of the aerial vehicle in the optimization problem. Since the ratio of communication energy to the total energy consumed by the aerial vehicle is negligible, the results proposed in the above works are suboptimal for the GEE maximization of ACS.  When the altitude of an AAP increases, the LoS coverage area increases, the LoS channel gain decreases and the energy consumed by the aerial vehicle also increases. With these facts, we can say that the GEE of an ACS, defined as the ratio of the total number of data bits transmitted to the total energy consumed, will not be maximum either at maximum or minimum permitted AAP altitudes. We exploit this tradeoff between the total number of bits transmitted and the energy consumed to determine an energy-efficient hovering altitude for the AAP. Some of the works which consider the aerial vehicle's energy consumption include \cite{eom2019uav}, \cite{rui2}, \cite{rui3}. An energy-efficient 3D trajectory of a UAV deployed to serve a set of IoT nodes is investigated in \cite{8038869}. Optimal trajectories, which minimize the fixed and rotary-wing UAV associated energy are designed in \cite{rui2} and \cite{rui3} respectively. The authors in \cite{eom2019uav} maximizes the minimum average rate and energy efficiency through joint optimization of trajectory, velocity, and acceleration of UAV flying at a fixed altitude. An altitude-dependent energy consumption model is used by the authors of \cite{ZorbasDimitrios2016Odpa} to find drone locations that minimize the cost while ensuring the surveillance of all the targets.

To the best of our knowledge, we are the first to determine an optimal altitude which maximizes the GEE for an ACS considering both the energy required for communication and energy consumed by the aerial vehicle. The rest of the paper is organized as follows. In section \ref{s1}, we model the system. The energy consumption of the aerial vehicle and GEE are explained in section \ref{geea}. Section \ref{s2} formulates the optimization problem and solves it using SCP and MFP techniques. The numerical results are discussed in section \ref{s3}. Finally, our findings are concluded in section \ref{s4}.

In this paper, scalars are represented by lowercase letters. Boldface lowercase letters are used to denote vectors and $\mathbb{R}^{M}$ denotes the set of $M$ dimensional real-valued vectors.
 
\section{System Model}
\label{s1}
\subsection{System Model}  

\begin{figure}
\centering
\includegraphics[width=0.9\linewidth]{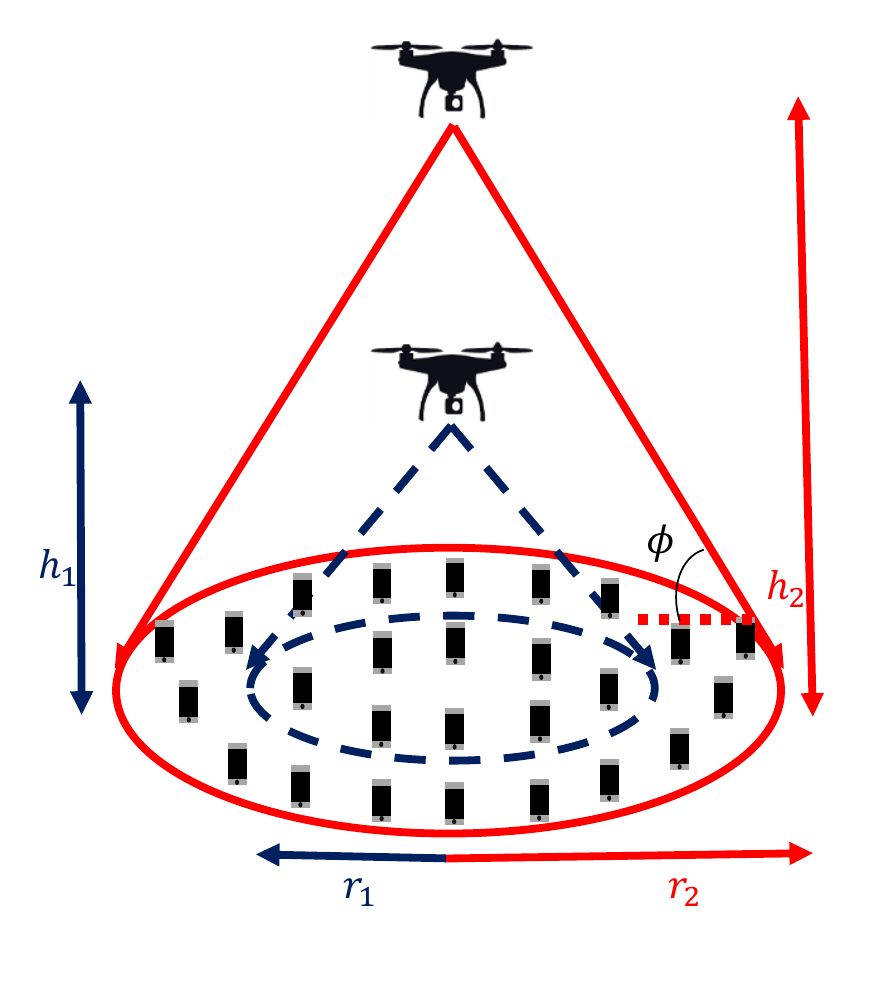}
\caption{AAP enabled downlink broadcast transmission scenario. }
\label{figg1}
\end{figure}

We consider an orthogonal multiple access downlink broadcast transmission scenario enabled by an AAP acting as a flying base station, where each user is allocated a fixed bandwidth. We assume there is always a sufficient number of orthogonal channels (e.g., narrowband frequency division multiple access systems \cite{goldsmith2005wireless}). As shown in Fig.\ref{figg1}, we assume a uniform distribution of $N$ UEs in the AAP coverage area $A_{ue}=\pi\overline{r}^{2}$ such that $N=\rho_{ue}\ast A_{ue}$, where $\rho_{ue}$ and 
$\overline{r}=h_{A}\text{cot}(\phi)$ represent the density of UEs and the radius of the AAP coverage area respectively, and $\phi$ represents the minimum elevation angle required for the LoS channel between the edge UE and the AAP \cite{angleforlos}. The AAP is employed at an altitude of $h_{A}$ meters (m) with the horizontal plane coordinates the same as the center of $A_{ue}$. In addition to this, we consider the deployment of this system in rural areas where the channel between the AAP and UE is dominated by the LoS link. In real life, this represents the access segment of an ACS in which an AAP is deployed for cellular coverage extension in a rural area. Given this, the LoS channel gain between the UE located at a distance $r$ from the center of the coverage area and the AAP 
is given by
\begin{IEEEeqnarray}{c}
h(r)=\dfrac{h_{0}}{r^{2}+h_{A}^{2}}
\end{IEEEeqnarray}
where $h_{0}$ represents the channel gain at a reference distance of 1m.
The signal-to-noise ratio (SNR) $\gamma_{\overline{r}}$, at the edge UE is given by
\begin{IEEEeqnarray}{c}
\gamma_{\overline{r}} = \dfrac{P_{T}h_{0}}{N(\overline{r}^{2}+h_{A}^{2})\sigma^{2}}
\label{snr}
\end{IEEEeqnarray}
where the total data transmission power $P_{T}$ available at the AAP is divided equally among the $N$ UEs and $\sigma^{2}$ represents the variance of the zero-mean additive white Gaussian noise at the corresponding receiver.

Hence, the total number of bits transmitted per unit Hz of bandwidth from the AAP to the considered UE through a channel of bandwidth $W$ in $T$ seconds is expressed as
\begin{IEEEeqnarray}{c}
R(\overline{r})=  T\text{log}_{2}(1+\gamma_{{\overline{r}}})
\label{R} \,\,\text{bits/Hz}
\end{IEEEeqnarray}
Through \eqref{R} and \eqref{snr},  the data rate of a UE  depends on distance $r$ from the center of the coverage area.
Because of the inverse relationship between $\gamma_{r}$ and $r$, the data rate of any UE is lower bounded by the data rate of the edge UE. That is
\begin{IEEEeqnarray}{c}
R(r)\geq R(\overline{r})=R(h_{A})\,\,\,\,\,\,\,\,\,\,\,\,\,\,\,   \forall r \,\,\leq\,\,\overline{r}
\label{Relation} 
\end{IEEEeqnarray}
Since maximizing $R(r)$ is equivalent to maximizing $R(\overline{r})$ and for ease of explanation, we consider the sum of minimum rate, $R(\overline{r})$, in the definition of the GEE of the considered ACS in section \ref{geea}. The algorithm developed in section \ref{s2} applies to the maximization of GEE defined in terms of the sum of actual rate $R(r)$.
\section{Global Energy Efficiency of the ACS}
\label{geea}
The global energy efficiency of the considered ACS is given by
 \begin{IEEEeqnarray}{c}
\text{GEE}[\text{bits}/{\text{Joule.Hz}}]=\dfrac{\overline{R}(h_{A})[\text{bits/Hz}]}{E(h_{A}, T)[\text{Joule}]}
\label{gee}
\end{IEEEeqnarray}
where $\overline{R}(h_{A})$ is the sum of the minimum number of data bits transmitted per Hz from the AAP to the $N$ UEs in $T$ seconds; $E(h_{A}, T)=E_{A}(h_{A}, T)+E_{C}(T)$ is the total energy consumed by the AAP, in which $E_{C}(T)$ is the energy required for data communication and $E_{A}(h_{A}, T)$ given by \eqref{Eaa}, is the total energy consumed by the mechanical parts of the aerial vehicle during vertical climb and hovering. We consider a climb-hover communicate scheme in which the AAP climbs at a specific altitude and then communicates with $N$ UEs while hovering.

\subsection{Sum of the minimum number of data bits transmitted,  $\overline{R}(h_{A})$ }
Considering the uniform distribution of UEs over $A_{ue}$, the sum of the minimum number of data bits transmitted per Hz from the AAP to the $N$ UEs in $T$ seconds through orthogonal channels of bandwidth $W$ Hz is expressed as
\begin{IEEEeqnarray}{rCl}
\overline{R}(h_{A})& = & T \int_{0}^{2\pi }\int_{0}^{\bar{r}}\rho_{ue}R(\overline{r}) r drd\theta\nonumber\\
%& = & T \rho_{ue}\int_{0}^{2\pi }\int_{0}^{h_{A}cot\phi}{\text{log}_{2}\left ( 1+\dfrac{P_{T}h_{0}}{N({\overline{r}}^{2}+h_{A}^{2})\sigma^{2}} \right )}r dr d\theta
%\nonumber\\
& = & T \rho_{ue} \pi h_{A}^{2} \text{cot}^{2}\phi{\text{log}_{2}\left ( 1+\dfrac{P_{T}h_{0}}{N({\overline{r}}^{2}+h_{A}^{2})\sigma^{2}} \right )}
\nonumber\\
& = & T \rho_{ue} \pi h_{A}^{2} \text{cot}^{2}\phi \text{log}_{2}\left (1+ \frac{\beta} {h_{A}^{4}}\right )
\nonumber
\label{sumrate}
\end{IEEEeqnarray}
where $\beta=\dfrac{P_{T}h_{0}\text{sin}^{2}\phi}{\pi\rho_{ue}\text{cot}^{2}\phi\sigma^{2}}$.
\subsection{Aerial vehicle energy consumption}
The total energy consumed by an aerial vehicle ($E(h_{A}, T)$) is composed of three main parts:
\begin{enumerate}
\item Energy required for data communication ($E_{C}(T)$).  
  \item Energy consumed by the rotor of the aerial vehicle during climbing from ground to an altitude of $h_{A}$ $(E_{cl}(h_{A}))$.
  \item Energy consumed by rotor during hovering at altitude $h_{A}$ $(E_{ho}(h_{A},T))$.  
\end{enumerate}
The energy required for data communication is given by 
\begin{IEEEeqnarray}{rCl}
E_{C}(T)=(P_{T}+P_{H})T
\end{IEEEeqnarray}
where $P_{T}$ is the total power used for the symbol transmission and $P_{H}$, is the total power consumption by all the hardware circuits in the transmitter section of the AAP.
 The energy parts $E_{cl}(h_{A})$ and $E_{ho}(h_{A},T)$ follow the energy consumption model presented by the authors of \cite{energyquad}. In \cite{energyquad},  the authors presented different power/energy consumption factors based on the field experiments performed on the Intel Aero Ready to Fly Drone. Unlike fixed and rotary-wing unmanned aerial vehicles \cite{rui2},\cite{rui3}, the energy consumed by the rotor of a quadropter/drone during hovering is dependent on the hovering altitude \cite{energyquad},\cite{ZorbasDimitrios2016Odpa}. According to \cite{energyquad}, the energy consumed by the quadropter during climbing from the ground to an altitude of $h_{A}$ with a constant climb rate is given by
\begin{IEEEeqnarray}{rCl}
E_{cl}(h_{A})= \alpha_{cl}h_{A}+\beta_{cl}
\end{IEEEeqnarray}
and the energy consumed during hovering at an altitude $h_{A}$ for $T$ seconds is given by
\begin{IEEEeqnarray}{rCl}
E_{ho}(h_{A},T)= (\alpha_{ho}h_{A}+\beta_{ho})T
\end{IEEEeqnarray}
where the constants $\alpha_{cl}, \beta_{cl}, \alpha_{ho}, \beta_{ho}$ are determined from the curve fitting performed on the measured power/energy values from the field experiments.

\begin{figure}
\centering
\includegraphics[width=\linewidth]{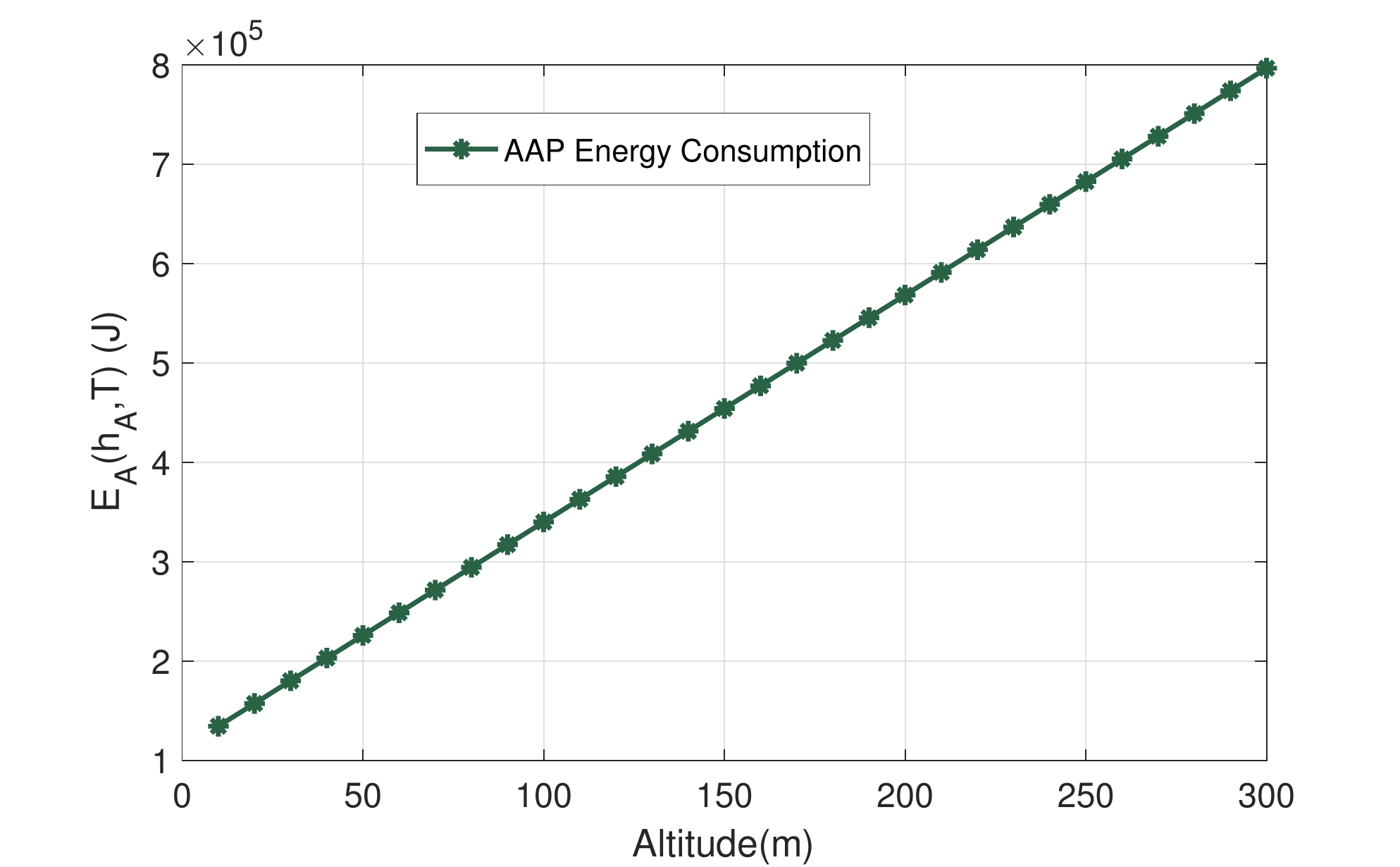}
\caption{Total energy consumed by the rotor of the aerial vehicle.}
\label{fig1}
\end{figure}

Hence the total energy consumed by the rotor of the AAP to climb to an altitude of $h_{A}$ m and hover for $T$ seconds is given by
\begin{IEEEeqnarray}{c}
E_{A}(h_{A}, T)=E_{cl}(h_{A})+E_{ho}(h_{A},T)
\label{Eaa}
\end{IEEEeqnarray}
Figure \ref{fig1} shows the increasing nature of $E_{A}(h_{A}, T)$ with altitude for a fixed time of operation with  constants   $\alpha_{cl}=315,\beta_{cl}=-211.261,\alpha_{ho}=4.917,\beta_{ho}=275.204$ \cite{energyquad} and $T=400s$. It is because as the altitude increases, the air temperature and pressure decreases. The decrease in the air pressure reduces the upward thrust provided by the air, to balance the downward force produced by the weight of the aerial vehicle. Hence, at higher hovering altitudes, to balance the weight, the propeller of the aerial vehicle needs to generate an additional force, which results in increased energy consumption.

Hence the total energy consumed by the AAP is given by
\begin{IEEEeqnarray}{c}
E(h_{A},T)=E_{cl}(h_{A})+E_{ho}(h_{A},T)+E_{C}(T)
\label{Ea}
\end{IEEEeqnarray}
\section{Problem Formulation}
\label{s2}
Our objective is to find the optimum altitude for the AAP, which maximizes the system's global energy efficiency (GEE) subject to minimum data rate and altitude constraints. By using \eqref{gee},\eqref{sumrate} and \eqref{Ea}, our main objective is formulated as an optimization problem and is expressed as follows:
 \begin{IEEEeqnarray}{rCl}
\text{(P1)} & : & \underset{h_{A}}{\text{maximize}}\,\,\,\,   \dfrac{\overline{R}(h_{A})}{E(h_{A}, T)}\nonumber\label{p1}\\
\text{s.t.} & & h_{min}\leq h_{A}\leq h_{max}\label{c1}\\
& & W \text{log}_{2}\left(1+\dfrac{\beta}{h_{A}^{2}(\overline{r}^{2}+h_{A}^{2})}\right) \geq  R_{0}\label{c2} 
\end{IEEEeqnarray}
where \eqref{c1} represents the permitted AAP altitude range specified by the aviation regulatory board and $R_{0}$ is the minimum data rate required by the UE in bits-per-second (bps). 

The objective function of (P1) belongs to the class of fractional programming problems \cite{fractional}. It can  be globally solved using Dinkelbach's algorithm \cite{DB},  provided $\overline{R}(h_{A})$ is concave, and $E(h_{A}, T)$, \eqref{c1}, \eqref{c2} are convex functions of $h_{A}$. From \eqref{Ea}, \eqref{c1}, we find that the denominator of the objective function and AAP altitude constraint are a convex function of $h_{A}$ while the numerator $\overline{R}(h_{A})$ in \eqref{sumrate}, is neither convex nor concave in nature. Also, the minimum individual data rate constraint \eqref{c2} is non-convex in nature. Hence (P1) cannot be globally solved with polynomial time complexity. As a means to obtain an efficient solution that fulfills the Karush Kuhn Tucker (KKT) conditions of (P1), we use the polynomial-time complex sequential convex programming (SCP) technique \cite{Marks}. Besides, to obtain the global optimum of (P1), we exploit the monotonic structure of the objective function in the general framework of monotonic fractional programming (MFP) optimization \cite{scp}\cite{monotonic}. 
\subsection{GEE Maximization using SCP}
In this part, we find the optimal altitude of the AAP, which maximizes the GEE of the ACS using sequential convex programming. The fundamental idea of SCP is to iteratively solve a sequence of convex approximated problems of the original non-convex problem so that the feasible solution points converge to the KKT point of the original non-convex problem \cite{Marks}. Here we approximate the non-concave numerator, $\overline{R}(h_{A})$ of (P1) as a concave function using a first-order Taylor approximation technique.

For the $k^{th}$ iteration, let $h_{k}$ be the feasible solution from the previous iteration. Then the first order Taylor approximation of $\overline{R}(h_{A})$ about $h_{k}$ is
\begin{IEEEeqnarray}{rCl}
\overline{R}(h_{A}) & \approx &\overline{R}(h_{k})+{\overline{R}}'(h_{k})(h_{A}-h_{k})
\label{taylor1}
\end{IEEEeqnarray}
where
\begin{IEEEeqnarray}{rCl}
{\overline{R}}'(h_{k}) & = &  T \rho_{ue} \pi \text{cot}^{2}\phi 2 h_{k} \text{log}_{2}\left (1+ \frac{\beta} {h_{k}^{4}}\right )\nonumber \\
&& - T \rho_{ue} \pi \text{cot}^{2} \phi \dfrac{4\beta h_{k} }{\text{log}_{e}(2)\left({\beta}+{h_{k}}^{4}\right)}\\
%\label{taylor2}
\end{IEEEeqnarray}
Similarly, the non-convex nature of \eqref{c2} is tackled with the following Taylor approximation:
\begin{IEEEeqnarray}{c}
W \text{log}_{2}\left(1+\dfrac{\beta}{h_{k}^{4}}\right)-\dfrac{4 W \beta (h_{A}-h_{k})}{ h_{k}(\beta+h_{k}^{4})\text{log}_{e}2} \geq R_{0}
\label{c3}
\end{IEEEeqnarray}
Using \eqref{taylor1} and \eqref{c3}, (P1) can be reformulated as 
\begin{IEEEeqnarray}{rCl}
\text{(P2)} & : & \underset{h_{A}}{\text{maximize}}\,\,\,\,   \dfrac{\overline{R}(h_{k})+{\overline{R}}'(h_{k})(h_{A}-h_{k})}{E(h_{A}, T)}\nonumber\\
\text{s.t.} & & \eqref{c1}, \eqref{c3}
\end{IEEEeqnarray}
\begin{algorithm}[]
\caption{GEE Maximization using SCP}
%\SetAlgoLined
Initialize $h_{1}$, $l_{1}^{s}=\dfrac{S(h_{1},h_{1})}{E(h_{1},T)}$, $k=1$.\\
 \While{$(1)$}
 {
 $h_{opt}^{s}=h_{k}$ \\
 Determine the optimal solution $h_{k}^{s*}$ by solving 
\begin{IEEEeqnarray}{rc}
&  \underset{h_{A}}{\text{maximize}}\,\,\,\, S(h_{A},h_{k})-l_{k}^{s}E(h_{A}, T)\nonumber\\
& s.t\,\,\, \eqref{c1}, \eqref{c3}\nonumber
\end{IEEEeqnarray} \\
$l_{k+1}^{s}=\dfrac{S(h_{k}^{s*},{h_{k}})}{E(h_{k}^{s*}, T)}$\\  
 \If{(($l_{k+1}^{s}-l_{k}^{s})/{l_{k+1}^{s}})<\zeta$}
 {
break;
 }
 
 $h_{k+1}=h_{k}^{s*}$\\
 $k=k+1$\\
 }
 \textbf{Output}:\,{Optimal AAP Altitude$ = h_{opt}^{s}$ }
 \label{algorithm1}
\end{algorithm}

Note that (P2) is a single ratio fractional maximization problem with a concave numerator $S(h_{A},h_{k})=\overline{R}(h_{k})+{\overline{R}}'(h_{k})(h_{A}-h_{k})$, convex denominator $E(h_{A},T)=E_{A}(h_{A}, T)+E_{C}(T)$ and convex constraints. Therefore (P2) can be efficiently solved by using polynomial time complex Algorithm 1. In every iteration of Algorithm 1, the optimal solution in step 4 is determined by using standard convex optimization tools like CVX \cite{cvx}. In section \ref{s3}, we show that the efficient solution of (P1) obtained by solving (P2) through Algorithm 1 matches the global optimum obtained using the monotonic fractional program optimization technique.
\subsection{GEE maximization using MFP}
The candidate solution obtained from SCP cannot be considered as the global optimum of (P1). Therefore, to obtain the global optimum of (P1), we exploit the monotonic behavior of the objective function using the monotonic fractional programming technique \cite{scp}\cite{monotonic}. The key idea is that the global optimum of an increasing objective function of a maximization problem lies in the outer boundary of the feasible set formed by the constraints.
Following the fundamental definitions from \cite{monotonic}, a maximization problem takes the canonical form of a monotonic optimization problem, if it can be formulated as
\begin{IEEEeqnarray}{rCl}
\text{(P3)} & : & \underset{\mathbf{v}}{\text{maximize}}\,\,\,\,   f(\mathbf{v})\nonumber\\
\text{s.t.} & & \mathbf{v} \in \mathcal{G} \cap \mathcal{H} \nonumber
\end{IEEEeqnarray}
where $f:\mathbb{R}^{M} \rightarrow \mathbb{R}$ is an increasing function of $\mathbf{v}$, $\mathcal{G}\subset[\mathbf{0},\mathbf{a}]$ is a
compact normal set with nonempty interior, and $\mathcal{H}$ is a closed conormal set on $[\mathbf{0},\mathbf{a}]$. For exact definitions of monotonicity, normal and co-normal sets please refer to \cite{monotonic}.

The optimization problem (P1) fits in the class of fractional problems, which can be globally solved by Algorithm 2. For a given positive $l_{k}^{m}$, in every $k^{th}$ iteration of Algorithm 2, we need to solve the following maximization problem in step $4$:
\begin{IEEEeqnarray}{rCl}
\text{(P4)}&:& \underset{h_{A}}{\text{maximize}}\,\,\,\,\overline{R}(h_{A})-l_{k}^{m}\left\{E_{A}(h_{A}, T)+E_{C}(T)\right\}\IEEEeqnarraynumspace\label{m1}\\
\text{s.t.} & & \eqref{c1} - \eqref{c2}
\end{IEEEeqnarray}
 It should be noted that, at first look, (P4) doesn't take the canonical form of monotonic optimization problem defined in (P3). However, (P4) can be expressed as the maximization of differences of increasing functions of $h_{A}$, which allows us to reformulate (P4) as a monotonic optimization problem. For the ease of reformulation, we equivalently represent the minimum individual data rate constraint as
 \begin{IEEEeqnarray}{c}
 h_{max}=\left[ \dfrac{\beta}{2^{\dfrac{R_{0}}{W}}-1}\right]^{1/4}
 \label{hmax} 
 \end{IEEEeqnarray}
  Note that \eqref{m1} can be rewritten as 
  \begin{IEEEeqnarray}{rCl}
 \text{(P5)} & : & \underset{h_{A}}{\text{maximize}}\,\,\,\,\overline{R}_{1}(h_{A})- \overline{R}_{2}(h_{A},l_{k}^{m})\label{m2}\\
  \text{s.t.} & & \eqref{c1}
\end{IEEEeqnarray}
where
\begin{IEEEeqnarray}{rCl}
 \overline{R}_{1}(h_{A}) & = & T \rho_{ue} \pi \text{cot}^{2}\phi h_{A}^{2} \text{log}_{2}\left (\beta +h_{A}^{4} \right )\nonumber\\
  \overline{R}_{2}(h_{A},l_{k}^{m}) & = & T \rho_{ue} \pi \text{cot}^{2}\phi h_{A}^{2} \text{log}_{2}\left (h_{A}^{4} \right )
\label{m3}\nonumber\\
& & +\> l_{k}^{m}(E(h_{A}, T)) 
\end{IEEEeqnarray}
are monotonically increasing functions of $h_{A}$, and $h_{max}$ of \eqref{c1} is given by \eqref{hmax}. In order to write (P5) in canonical form, we introduce the additional variable $t=\overline{R}_{2}(h_{max},l_{k}^{m})-\overline{R}_{2}(h_{A},l_{k}^{m})$, which allows (P5) to be reformulated as
 \begin{IEEEeqnarray}{rCl}
\text{(P6)}  & : & \underset{h_{A},t}{\text{maximize}}\,\,\,\,\overline{R}_{1}(h_{A})+t\label{m2}\\
\text{s.t.} & & (h_{A}, t)\in \mathcal{G}\cap \mathcal{H}
\end{IEEEeqnarray}
where
 \begin{IEEEeqnarray}{l}
 \mathcal{G} =\left\{ 
\begin{aligned} 
 (h_{A}, t):h_{A}\leq h_{max},\\
   t\leq \overline{R}_{2}(h_{max},l_{k}^{m})-\overline{R}_{2}(h_{A},l_{k}^{m}) \label{g} 
\end{aligned} \right\} \\
\nonumber\\
\nonumber\\
\mathcal{H} =  \left\lbrace (h_{A}, t):h_{A}\geq h_{min}, t \geq 0\right\rbrace \label{h}
 \end{IEEEeqnarray}
 By the monotonically increasing behavior of $\overline{R}_{2}(h_{A},l_{k}^{m})$ we can relate
 \begin{IEEEeqnarray}{r}
 \overline{R}_{2}(h_{min},l_{k}^{m})\leq\overline{R}_{2}(h_{A},l_{k}^{m})
\end{IEEEeqnarray} 
 By [proposition 2, \cite{scp}], \eqref{g} defines a normal set and \eqref{h} defines a co-normal set in the polyblock
 \begin{IEEEeqnarray}{r}
 \left[ h_{min},  h_{max}\right] \times \left[0, \overline{R}_{2}(h_{max},l_{k}^{m})-\overline{R}_{2}(h_{min},l_{k}^{m})  \right]
\IEEEeqnarraynumspace
\label{poly}
\end{IEEEeqnarray}
with the vertex set $\mathcal{V}$. Hence by using \eqref{m2}-\eqref{h} we represent (P4) in the canonical form of monotonic optimization problem with $f(\mathbf{v})=\overline{R}_{1}(\mathbf{v}(1))+\mathbf{v}(2)$, $\forall\mathbf{v}\in\mathcal{V}$,  which can be globally solved by using the polyblock outer approximation algorithm as explained in Algorithm 3 \cite{monotonic}. Even though the complexity of this global optimization algorithm is exponential in the number of variables, it is much lower compared to other global optimization techniques, which exhaustively search over the entire feasible set. Hence the globally optimal AAP altitude is obtained by solving (P1) using Algorithm 2 in which, at each iteration, step 4 is solved by using Algorithm 3. 
\begin{algorithm}[]
\caption{GEE Maximization using MFP}
%\SetAlgoLined
Initialize $h_{1}$, $l_{1}^{m}=\dfrac{\overline{R}(h_{1})}{E(h_{1}, T)}$, $k=1$.\\
 \While{$(1)$}
 {
 $h_{opt}^{m}=h_{k}$ \\
 Determine the optimal solution $h_{k}^{m*}$ by solving the monotonic optimization problem (P6) using Algorithm \ref{PA}\\
$l_{k+1}^{m}=\dfrac{\overline{R}(h_{k}^{m*})}{E(h_{k}^{m*},T)}$\\  
 \If{(($l_{k+1}^{m}-l_{k}^{m})/{l_{k+1}^{m}})<\zeta$}
 {
break;
 }
 
 $h_{k+1}=h_{k}^{m*}$\\
 $k=k+1$\\
 }
 \textbf{Output}:\,{Optimal AAP Altitude$ = h_{opt}^{m}$} 
 \label{algorithm2}
\end{algorithm}
\section{Numerical Evaluation}
\label{s3}
In this section, we compare the optimal altitude values obtained through SCP and MFP optimization techniques. Furthermore, the convergence behavior of the PA algorithm; the impact of aerial vehicle's energy consumption on GEE; the variation of GEE with minimum data rate requirement are discussed. We consider $h_{0}=1.42\times10^{-4}$, $P_{T}=10$ dBm, $W=20$MHz, $\sigma_{0}^{2}=-169 \text{dBm}/\text{Hz}$ , $\phi=43^{\circ}$, $\rho=0.005 \text{UEs}/\text{m}^{2}$, $P_{H}=5$W, $T=400$s, $R_{o}=20$Mbps and $h_{min}=10$m.

\begin{figure}
\centering
\includegraphics[width=\linewidth]{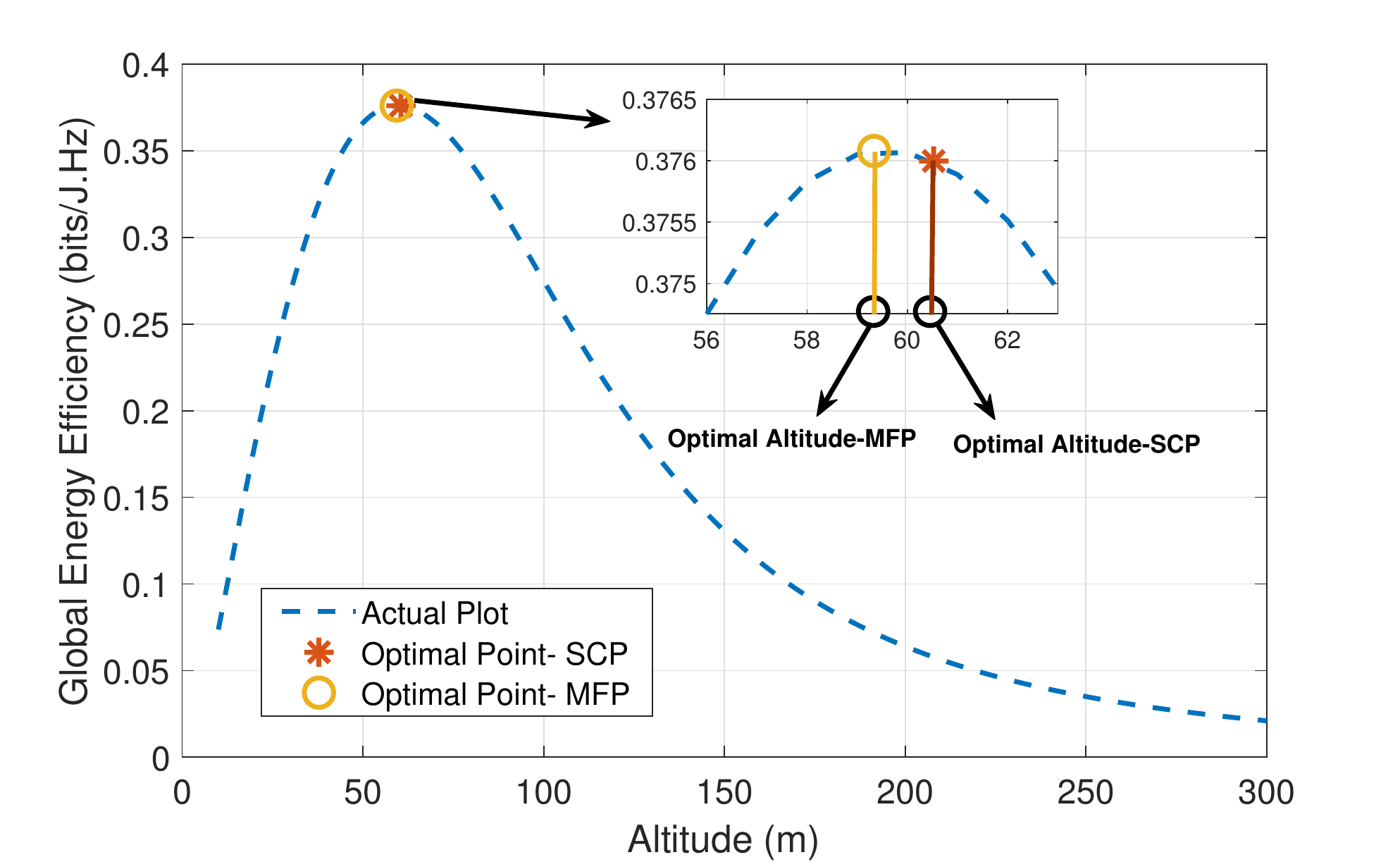}
\caption{Optimal solutions obtained from SCP and MFP. }
\label{fig3}
\end{figure}

Figure \ref{fig3} shows the accurate plot of GEE with the altitude of aerial vehicle, $h_{A}$, along with the optimal points obtained through SCP and MFP techniques. From the plot, it is observed that GEE decreases when $h_{A}$ is very low or very high. The reason for this behavior is that, at low $h_{A}$, the number of UEs covered ($N=\rho_{ue}\pi h_{A}^{2}\text{cot}^{2}\phi$) by the AAP decreases with decreasing $h_{A}$, leading to a decrease in the total number of bits transmitted, thereby to reduced GEE. At high altitude regions, the LoS channel gain between the UE and AAP decreases, the number of UEs covered by AAP increases, and $E(h_{A}, T)$ increases. In addition to this, with an increase in the number of users, power alloted for a single UE decreases. So in the high altitude region, the increase in the number of UEs is highly compensated by the combined effect of the decrease in channel gain, decrease in power per UE and increase in $E(h_{A}, T)$, which result in a low GEE. Figure  \ref{fig3} also shows that the optimal AAP altitude obtained by the SCP is very close to the globally optimal altitude obtained from the monotonic fractional programming technique. Hence the global optimum of our objective can be obtained by the polynomial-time complex sequential convex optimization technique. 
\begin{figure}
\centering
\includegraphics[width=\linewidth]{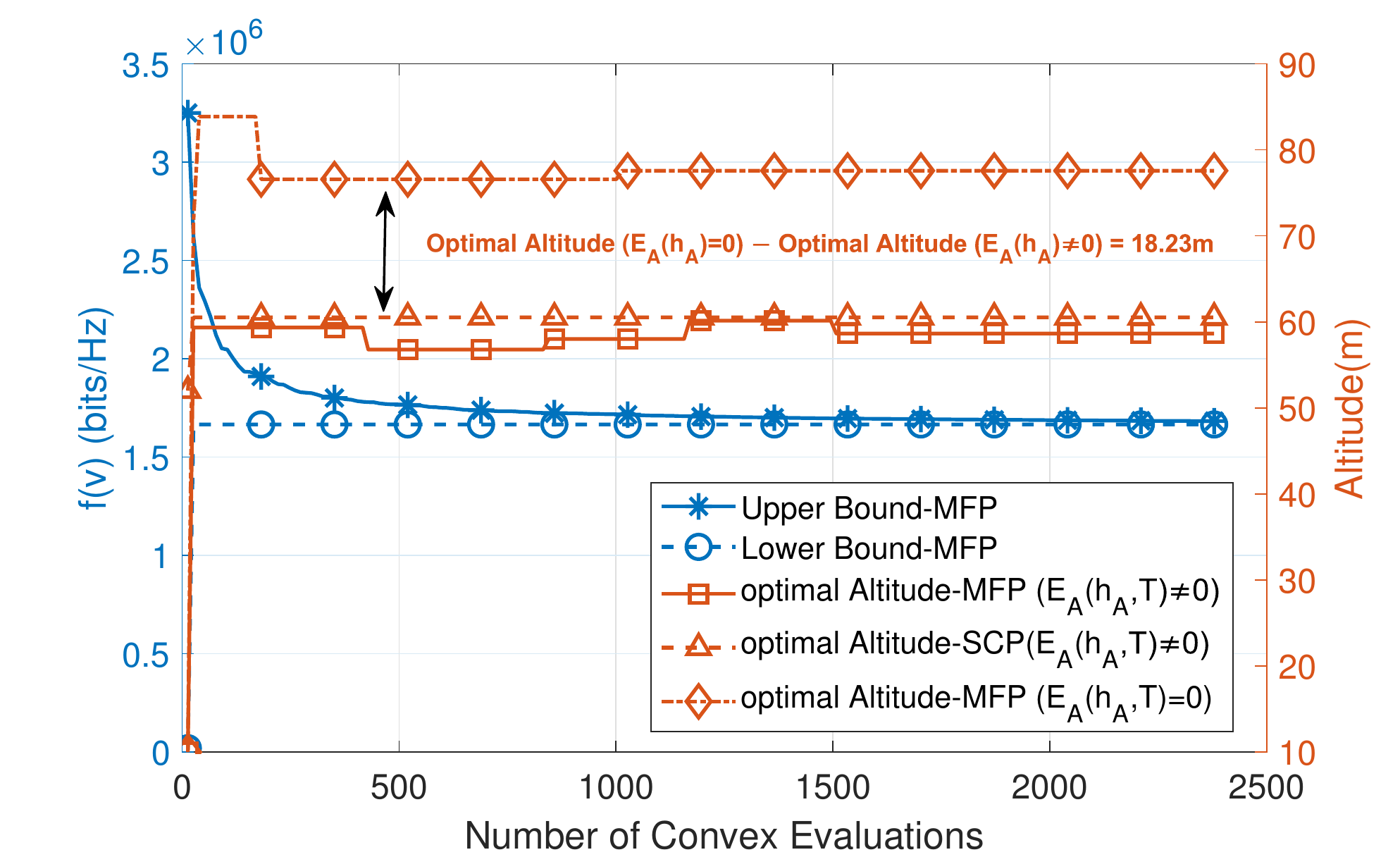}
\caption{Convergence behavior of the PA algorithm.}
\label{fig5}
\end{figure}

Figure \ref{fig5} shows the convergence behavior of the PA algorithm in the last iteration of MFP. The converging nature of upper($f_{max}$) and lower($f_{min}$) bounds of the PA algorithm guarantees the evaluation of the global optimum of GEE in a finite number of convex evaluations; with the number of convex evaluations much greater than that required by SCP. The optimal altitude plots of MFP and SCP says that the locally optimal altitude value obtained using SCP is equal to the globally optimal altitude obtained using MFP. In addition to this, Figure \ref{fig5} shows the error in determining the optimal altitude without considering the rotor energy consumption, $E_{A}(h_{A},T)$. It is observed that the optimal altitude determined with $E_{A}(h_{A}, T)=0,$ is 18.23m higher than the actual optimal altitude value. Hence, according to Figure \ref{fig3}, hovering at an altitude higher than the actual optimal value yields low GEE. Therefore, to achieve the maximum GEE value, the rotor energy consumption of the aerial vehicle should be considered while formulating the optimization problem. 
\begin{figure}
\centering
\includegraphics[width=\linewidth]{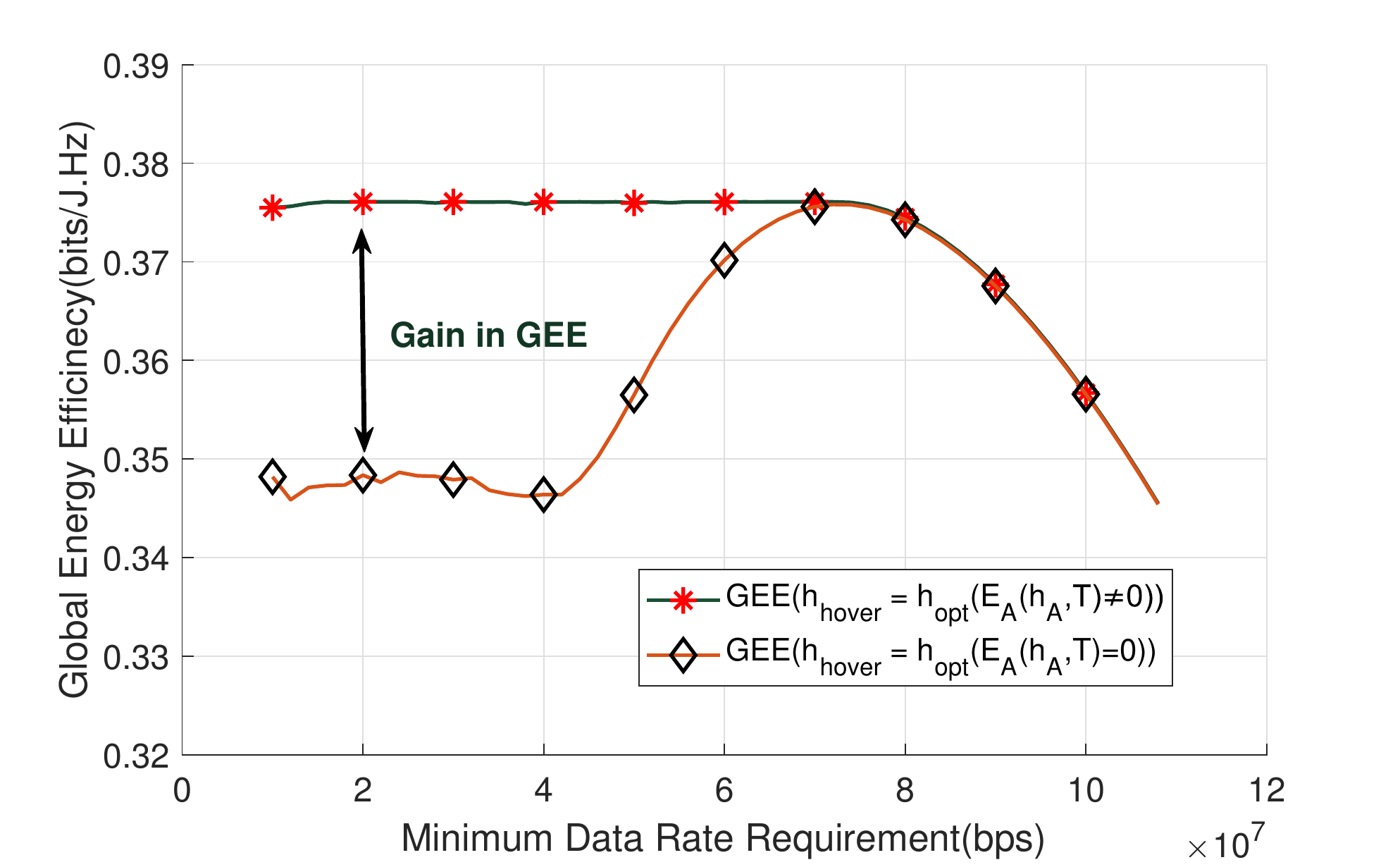}
\caption{Global energy efficiency versus minimum individual data rate requirement.}
\label{fig4}
\end{figure}

Figure \ref{fig4} depicts the variation of GEE with the minimum data rate requirement. The two GEE plots correspond to the variation of GEE with $R_{o}$ when the aerial vehicle is hovering at the optimal altitude obtained: (a) with zero energy consumed by the rotor of the aerial vehicle($h_{hover}=h_{opt}(E_{A}(h_{A}, T)=0)$ ); (b) with non-zero energy consumed by the rotor of the aerial vehicle ($h_{hover}=h_{opt}(E_{A}(h_{A}, T)\neq0)$ ). As said before, the GEE of the ACS with $h_{hover}=h_{opt}(E_{A}(h_{A}, T)\neq0)$ is more than that with $h_{hover}=h_{opt}(E_{A}(h_{A}, T)=0)$. This gain in GEE reflects the effect of considering the aerial vehicle's rotor energy consumption in altitude optimization, showcasing the novel aspect of this paper. In the plot, the value of GEE is constant for a range of $R_{o}$ and then it starts decreasing with an increase in $R_{o}$. From \eqref{hmax}, the value of the maximum allowed altitude, $h_{max}$ decreases with increase in $R_{o}$. It is because a higher minimum individual data rate is achieved by increasing the LoS channel gain obtained by decreasing the hovering altitude of the aerial vehicle. When $h_{max}(R_{o})$ is greater than the $h_{A}$ corresponding to the global optimum of GEE ($\text{GEE}_{\text{global}}$), the optimal altitude is equal to $h_{A}(\text{GEE}_{\text{global}})$ and GEE remains constant; for $h_{max}(R_{o})\leq h_{A}(\text{GEE}_{\text{global}})$, optimal altitude is equal to $h_{max}$, results in decrease in GEE with increase in $R_{o}$. The decrease in GEE with increase in $R_{o}$(decrease in $h_{max}$) shows the monotonically increasing property of GEE which is exploited in MFP.
\section{Conclusion} 
\label{s4}
In this work, we found the optimal energy-efficient altitude of an aerial access point which acts as a flying base station for an orthogonal multiple access downlink broadcast transmission scenario. The modeled energy consumption is the sum of energy consumed by the aerial vehicle and the energy required for the communication between the AAP and the UEs. An efficient solution to the formulated GEE maximization problem with individual data rate constraint and altitude constraint is obtained using sequential convex programming and is compared to the global optimum achieved by the monotonic fractional programming technique. One can see that the optimal altitude value from the polynomial-time complex SCP matches the globally optimal altitude value obtained from the monotonic fractional programming. Further, we observed that there is a gain in the GEE when the aerial access point is hovering at an optimal altitude determined by considering the non-zero rotor energy consumption of the aerial vehicle. In addition to this, the optimal altitude, and hence GEE, decrease with an increase in the minimum individual data rate constraint. Joint altitude and power optimization in a non-orthogonal multiple access transmission scheme with multiple AAPs is left as our future work.

\begin{algorithm}[]
\DontPrintSemicolon
  
%  \KwData{Testing set $x$}
Initialize $i=1$, $\mathcal{V}_{i}$ as the vertexset of polyblock \eqref{poly}\\
   Set $\mathbf{v}_{min}=argmin \lbrace f(\mathbf{v})\mid \mathbf{v}\in\mathcal{V}_{i} \rbrace$ \\
Set $\mathbf{v}_{max}=argmax \lbrace f(\mathbf{v})\mid \mathbf{v}\in\mathcal{V}_{i} \rbrace$\\
Set $f_{max}=max_{\mathbf{v}\in\mathcal{V}_{i}}f(\mathbf{v})$ and 
$f_{min}=f(\mathbf{v}_{min})$\\
\While {$((f_{max}-f_{min})/f_{max}>e)$} 
{
Obtain $\mathbf{v}_{o}$, the intersecting point of the line drawn from $\mathbf{v}_{min}$ to $\mathbf{v}_{max}$ with the normal region $\mathcal{G}$ using bisection method [Algorithm 1 \cite{bjornson2013optimal}]\\
Update the vertex set, $\mathcal{V}_{i+1}$ according to Lemma 2.16 of \cite{bjornson2013optimal}\\
\If {$f(\mathbf{v}_{o})>f_{min})$}
{
$f_{min}=f(\mathbf{v}_{o})$\\
$\mathbf{v}_{min}=\mathbf{v}_{o}$
}
set $i=i+1$\\
remove all $\mathbf{v}\in \mathcal{V}_{i}$ with $ f(\mathbf{v})\leq f_{min} +e $\\ 

Set $f_{max}=max_{\mathbf{v}\in\mathcal{V}_{i}}f(\mathbf{v})$\\

}
\textbf{Output}: $h_{k}^{m*}=\mathbf{v}_{o}(1)$.\\

%  \tcp*{this is a comment}
%  \tcc{Now this is an if...else conditional loop}
%  \If{Condition 1}
%    {
%        Do something    \tcp*{this is another comment}
%        \If{sub-Condition}
%        {Do a lot}
%    }
%    \ElseIf{Condition 2}
%    {
%    	Do Otherwise \;
%        \tcc{Now this is a for loop}
%        \For{sequence}    
%        { 
%        	loop instructions
%        }
%    }
%    \Else
%    {
%    	Do the rest
%    }
%    
%    \tcc{Now this is a While loop}
%   \While{Condition}
%   {
%   		Do something\;
%   }

\caption{PA Algorithm \cite{bjornson2013optimal}}
\label{PA}
\end{algorithm}
\thanks{Manuscript submitted April 2020. This work is supported by the project PAINLESS which has received funding from the European Union’s Horizon 2020 research and innovation programme under grant agreement No 812991.}
%\thanks{N. Babu and C.B.Papadias are with Research, Technology and Innovation Network (RTIN), 
%The American College of Greece, Greece (e-mail: nbabu@acg.edu, cpapadias@acg.edu).}
%\thanks{N. Babu, C.B.Papadias and P.Popovski are with Department of Electronic Systems, Aalborg University, Denmark (e-mail: niba,cop,petarp@es.aau.dk).}

\bibliographystyle{IEEEtran}
\bibliography{IEEEabrv,ff}
\end{document}